%% file: PolytopesFromGraphs.tex
\documentclass[12pt]{article}
\usepackage{amssymb}
\usepackage{amsmath}
\usepackage{graphics}
\usepackage{epsfig}
\usepackage{boxedminipage}

\newtheorem{thm}{Theorem}[section]
\newtheorem{lem}[thm]{Lemma}
\newtheorem{prp}[thm]{Proposition}
\newtheorem{cor}[thm]{Corollary}

\newtheorem{con}[thm]{Conjecture}
\newenvironment{proof}{{\it Proof}.\ }{\hfill$\square$\par\medskip}

\def\cF{{\mathcal F}}

\def\cR{{\mathcal R}}

\newcommand\SetOf[2]{\bigl\{#1\,\bigm|\,#2\bigr\}}

\def\dual{{\operatorname{dual}}}

\def\SetOf#1#2{\left\{\left.#1\vphantom{#2}\ \right|\ #2\vphantom{#1}\right\}}

\title{Reconstructing a non-simple polytope from its graph}
\author{Michael Joswig}

\begin{document}
\maketitle

\begin{abstract}
  A well-known theorem of Blind and Mani~\cite{Blind87:287} says that every simple polytope is
  uniquely determined by its graph.  In~\cite{Kalai88:381} Kalai gave a very short and elegant proof
  of this result using the concept of acyclic orientations.  As it turns out, Kalai's proof can be
  suitably generalized without much effort.  We apply our results to a special class of cubical
  polytopes.


\end{abstract}

\section{Introduction}

Polytopes arise as sets of admissible solutions for linear optimization problems.  Most practical
applications rely on Dantzig's Simplex Method to solve a given linear program.  This algorithm
starts at any vertex of the polytope and walks along improving edges to the optimum vertex.  Thus
the structure of the vertex-edge-graph is of considerable importance for the performance of the
Simplex Method.  There is reason to believe that geometric insight into the relationship between a
polytope and its graph will help to solve some of the many unanswered questions about the worst-case
running time of the Simplex Algorithm.  For an overview over facts about graphs of polytopes see
Ziegler~\cite[Chapter~3]{Ziegler98:0}.

In this context it is natural to ask how far the graph of a polytope determines the (combinatorial)
structure of the polytope.  The known answer is: not very much, in general.  For instance, there is
a class of simplicial $d$-polytopes with $n$~vertices, called {\em cyclic polytopes}, whose
vertex-edge-graph is the complete graph on $n$~vertices, which is the same as the graph of an
$(n-1)$-simplex.  In particular, in general, it is not even possible to determine the dimension of a
polytope from its graph. See also Gr\"unbaum~\cite[Chapter~12]{163.16603}.

The situation is quite different for simple polytopes, as it has been shown by Blind and
Mani~\cite{Blind87:287}.  Simple polytopes with isomorphic graphs are combinatorially equivalent.
The purpose of this paper is to show how much information in addition to the graph is necessary in
order to be able to reconstruct a general polytope.  The first part of the paper closely follows
Kalai's proof~\cite{Kalai88:381} of the theorem of Blind and Mani mentioned above.  For questions
concerning the complexity of the reconstruction procedure based on Kalai's proof see Achatz and
Kleinschmidt~\cite{Achatz99:0}.

The second part of this paper is devoted to a special class of cubical polytopes.  It is shown that
{\em capped\/} polytopes in the sense of Blind and Blind~\cite{970.45465} can be reconstructed from
their dual graphs.

Recently, Cordovil, Fukuda, and Guedes de Oliveira~\cite{CordovilFukudaOliveira} proved a Blind-Mani
type theorem in the context of oriented matroids.  And we should also mention a result of Bj\"orner,
Edelman, and Ziegler~\cite{698.51010}: A zonotope can be reconstructed from its graph.  The proofs
in both papers rely on entirely different techniques than the ones used here.

\section{The general reconstruction scheme}

Let $P$ be an arbitrary polytope and $\Gamma_P$ its vertex-edge graph.  For each vertex~$v$ of~$P$,
let $\eta_v$ be the map from the set~$\cF(v)$ of facets containing~$v$ to the set of subsets of the
set~$\Gamma_P(v)$ of edges through~$v$ such that each facet~$F$ is mapped onto the set of edges
through~$v$ which are contained in~$F$.  We call the family of images
$\SetOf{\eta_v(\cF(v))}{\mbox{$v$ vertex of~$P$}}$ the {\em edge labeled vertex figures\/}
of~$P$.  Observe that the vertex~$v$ can be recovered from the set~$\eta_v(\cF(v))$ of edge sets.

\begin{lem}\label{lem:vertex-figure}
  The face lattice of each vertex figure is determined by the edge labeled vertex figures.
\end{lem}

\begin{figure}[htbp]
  \begin{center}
    \epsfig{file=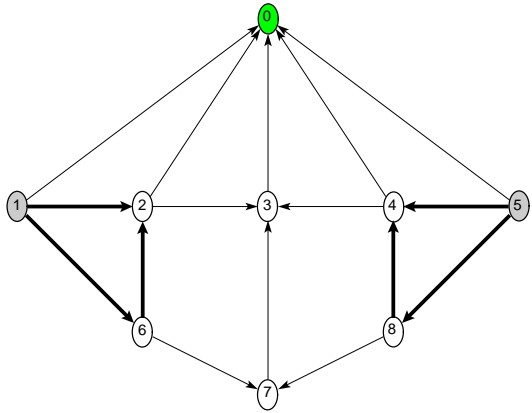,width=7cm}
    \caption[Non-AOF.]{A good acyclic orientation of the graph of a $3$-polytope which is not an
      abstract objective function.  The vertex~$0$ is the unique sink, while the vertices~$1$
      and~$5$ are sources.  The subgraph corresponding to the facets $\{1,2,6\}$ and $\{4,5,8\}$,
      respectively, is initial with respect to the acyclic orientation. }
    \label{fig:non-aof}
  \end{center}
\end{figure}

We follow Kalai's approach by examining acyclic orientations of~$\Gamma_P$, see also
Ziegler~\cite[Section~3.4]{Ziegler98:0}.  A {\em good\/} acyclic orientation has precisely one sink
on each non-empty face.  An {\em abstract objective function\/} of~$P$ is a good acyclic orientation
of~$\Gamma_P$ which also has precisely one source on each non-empty face.  For simple polytopes it
is known that the good acyclic orientations are precisely the {\em abstract objective functions\/}
and also the {\em shellings\/} of the boundary of the dual (simplicial) polytope.  This is a
consequence of Kalai's proof.  Note, however, that the graph of a non-simple polytope may have an
acyclic orientation which is {\em not\/} an abstract objective function; see the example in
Figure~\ref{fig:non-aof}.

An induced subgraph~$\Sigma$ is called {\em initial\/} with respect to an acyclic orientation if
there is no directed edge pointing to a vertex in~$\Sigma$ from any vertex of the complement
of~$\Sigma$.  For an example of an initial subgraph with respect to a good acyclic orientation
consider the following.  Take an arbitrary linear objective function with minimal face~$M$.  Tilting
the linear objective function slightly into general position induces a good acyclic orientation (and
even an abstract objective function).  The induced subgraph on~$M$ is initial.  In particular, good
acyclic orientations always exist.

In order to determine the combinatorial structure of~$P$ we have to find out which subsets of the
vertex set form facets.  Call a non-empty induced subgraph~$\Phi$ of~$\Gamma_P$ an {\em F-subgraph\/} if
it has the following properties:
\begin{enumerate}
\item There is a good acyclic orientation~$O$ of~$\Gamma_P$ such that $\Phi$ is initial with respect
  to~$O$.
\item For each vertex~$v$ of~$\Phi$ there is a (unique) facet $F\in\cF(v)$ such that
  $\eta_v(F)=\Phi(v)$, where by~$\Phi(v)$ we denote the set of edges through~$v$ in~$\Phi$.
\item The subgraph~$\Phi$ is minimal with respect to inclusion among the induced subgraphs of~$\Gamma_P$
  satisfying the properties above.
\end{enumerate}
Note that the minimality condition enforces the connectedness of an F-subgraph.  See again the
example in Figure~\ref{fig:non-aof} for an induced subgraph of the graph of a $3$-polytope
satisfying (i) and (ii), but not~(iii).

\begin{lem}\label{lem:F-subgraphs}
  The F-subgraphs of~$\Gamma_P$ are precisely the subgraphs of facets of~$P$.
\end{lem}

\begin{proof}
  Clearly, the induced subgraph on the vertex set of any facet is an F-subgraph.
  
  For the converse let $\Phi$ be an F-subgraph.  There is a good acyclic orientation~$O$
  of~$\Gamma_P$ such that $\Phi$ is initial with respect to~$O$. Choose a sink~$s$ among the
  vertices of~$\Phi$.  By assumption there is a facet $F\in\cF(s)$ with $\eta_s(F)=\Phi(s)$.
  Because $O$ is good, the vertex~$s$ is the unique sink of~$F$.  Moreover, all the vertices of~$F$
  are contained in~$\Phi$ because~$\Phi$ is initial with respect to~$O$.  Comparing $\Phi$ with the
  F-graph~$\Gamma_F\le\Gamma_P$ corresponding to~$F$ yields $\Phi=\Gamma_F$ due to the minimality
  of~$\Phi$.
\end{proof}

The proof of the preceding lemma shows that it is not necessary to require the uniqueness of the
facet in~(ii) of the definition of an F-subgraph.

\begin{thm}\label{thm:MainTheorem}
  A polytope can be reconstructed from its graph and the edge labeled its vertex figures.
\end{thm}

\begin{proof}
  Consider all possible acyclic orientations of the graph~$\Gamma_P$.  In view of
  Lemma~\ref{lem:F-subgraphs} it suffices to exhibit all good acyclic orientations of~$\Gamma_P$.
  
  Fix an acyclic orientation~$O$ and a vertex~$v$.  We want to compute the number~$f^O(v)$ of faces
  in which the vertex~$v$ is a sink with respect to~$O$.  This is the number of faces containing~$v$
  built from edges incoming at~$v$ only.  All faces through~$v$ can be enumerated because the face
  lattice of the vertex figure is known, cf.\ Lemma~\ref{lem:vertex-figure}.  Filter these faces for
  the incoming edge condition.  Note that for $P$ simple the number~$f^O(v)$ solely depends on the
  in-degree of~$v$ with respect to~$O$, which is not true for general polytopes.  In particular,
  we do not have a formula in closed form.
  
  Let $f^O=\sum_vf^O(v)$ and let $f$ be the number of non-empty faces. We proceed as Kalai in his
  proof.  As each non-empty face has at least one sink we have $f^O\ge f$.  Further, $O$ is good if
  and only if $f^O=f$.  Because good acyclic orientations always exist we have that
  $f=\min\SetOf{f^O}{\mbox{$O$ acyclic orientation on~$\Gamma_P$}}$.
\end{proof}

For $P$ a simple $d$-polytope, the edge labeled the vertex figure of the vertex~$v$ coincides with
the set $\binom{\Gamma_P(v)}{d-1}$, that is, the set of subsets of~$\Gamma_P(v)$ of
cardinality~$d-1$.  By iterated intersection they generate the boolean lattice~$[2]^d$ (top element
removed).  Each vertex figure of a simple $d$-polytope is a $(d-1)$-simplex.

\begin{cor}
  Every simple polytope can be reconstructed from its graph.
\end{cor}

\section{Cubical polytopes}

Simple polytopes are precisely the duals of simplicial polytopes.  From a certain perspective cubical
polytopes behave somewhat similar to simplicial polytopes.  So the following conjecture is tempting; 
see also~\cite[Problem~3.3]{CordovilFukudaOliveira}.

\begin{con}\label{Con}
  Every cubical polytope can be reconstructed from its dual graph.
\end{con}

The {\em dual graph\/}~$\Delta_P$ of a polytope~$P$ is the vertex-edge graph of the dual polytope
of~$P$, i.e.\ $\Delta_P=\Gamma_{P^\dual}$.  The nodes and edges of~$\Delta_P$ correspond to facets
and ridges (codimension~$2$ faces) of~$P$.

Note that, like for simplicial polytopes, there is no hope for cubical polytopes to be
reconstructible from their graph: In~\cite{NeighborlyCubical} it is shown that there are cubical
$d$-polytopes with the graph of the $n$-dimensional cube for arbitrary $n>d$.  The {\em neighborly 
  cubical polytopes\/} constructed in loc. cit. can be seen as cubical analogues of the (simplicial) 
cyclic polytopes mentioned in the introduction.

\begin{figure}[htbp]
  \begin{center}
   \begin{minipage}[c]{6cm}
     \epsfig{file=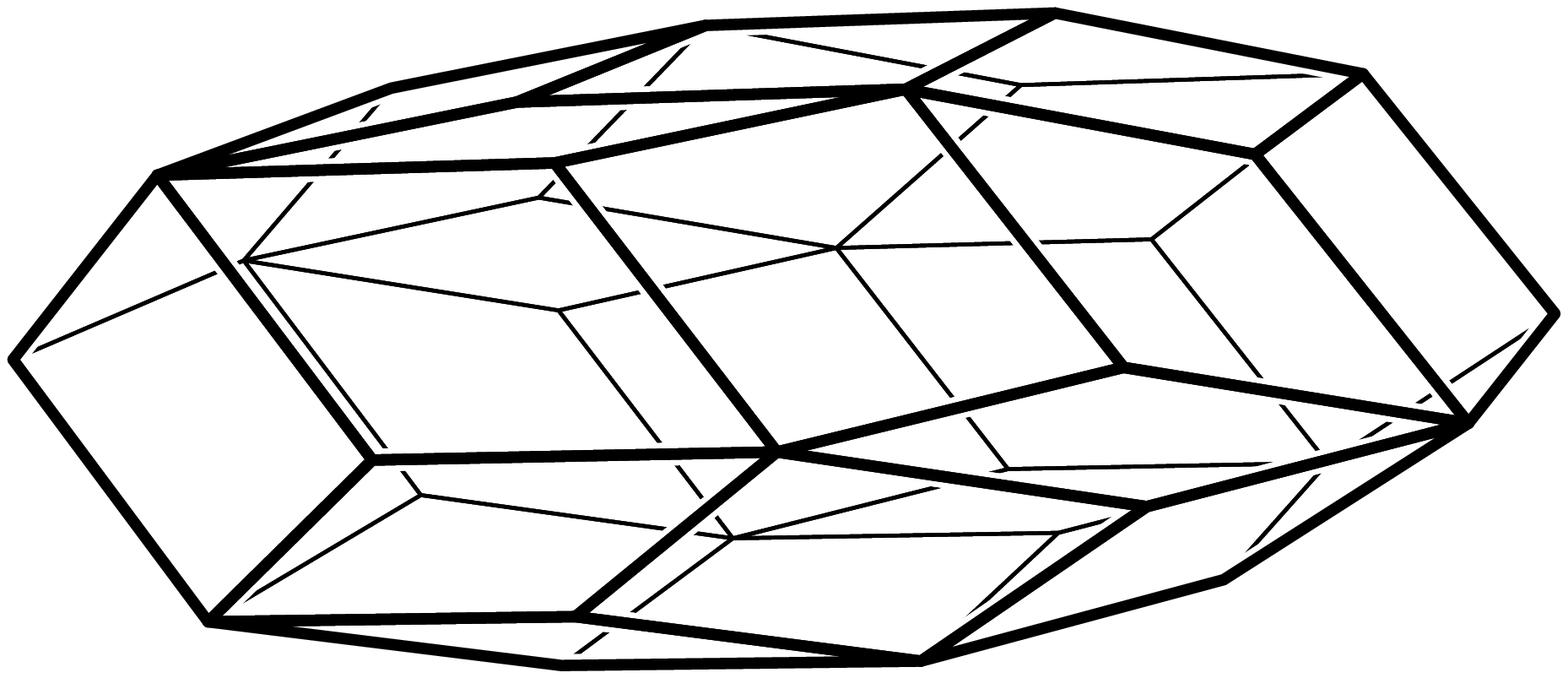,width=6cm}
    \end{minipage}\qquad
    \begin{minipage}[c]{6cm}
      \epsfig{file=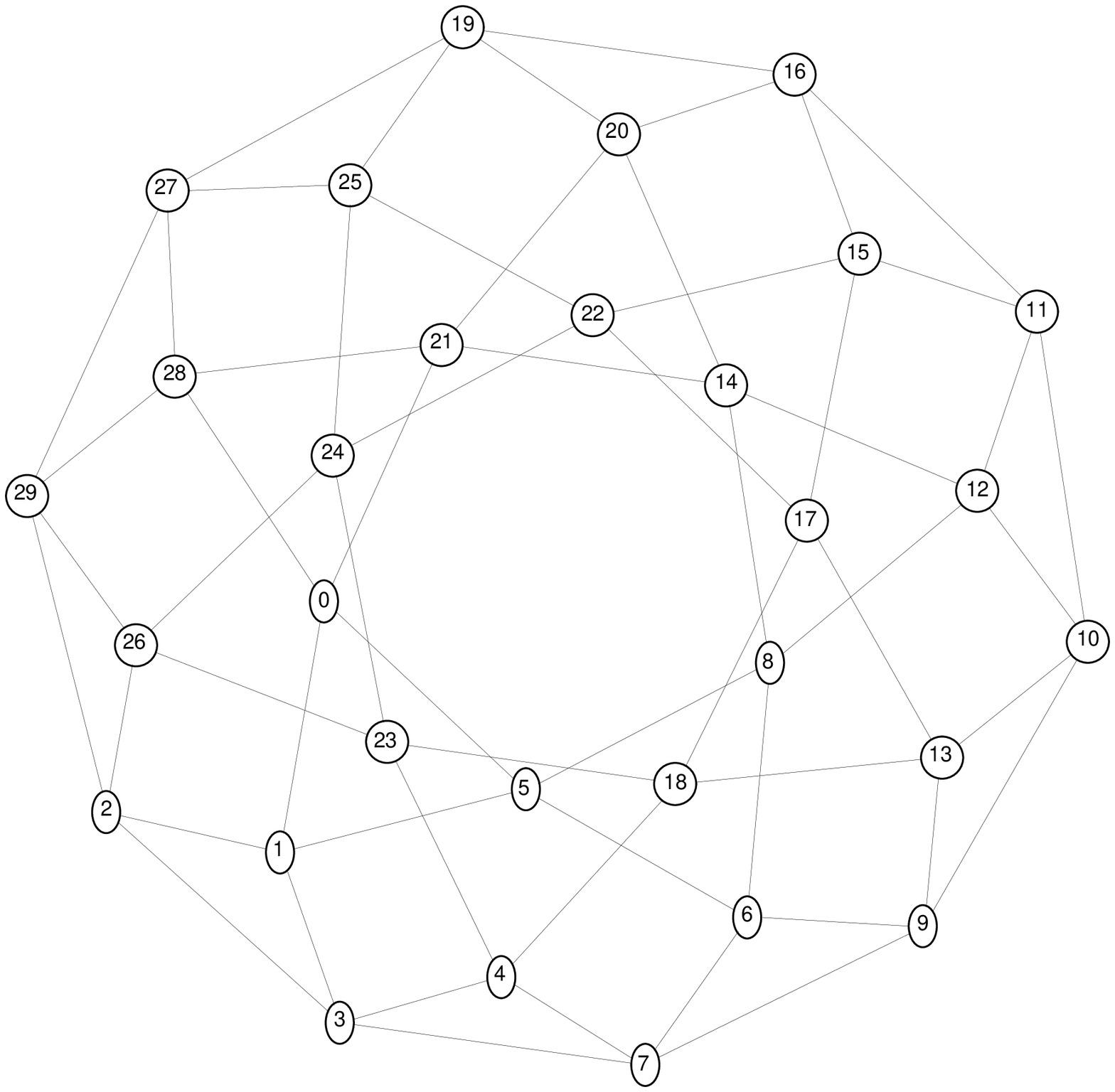,width=6cm}
    \end{minipage}
    \caption[Zonotope.]{A $3$-dimensional cubical zonotope (with $6$ zones) and its dual graph.}
    \label{fig:zonotope}
  \end{center}
\end{figure}

\begin{figure}[htbp]
  \begin{center}
    \begin{minipage}[c]{6cm}
      \epsfig{file=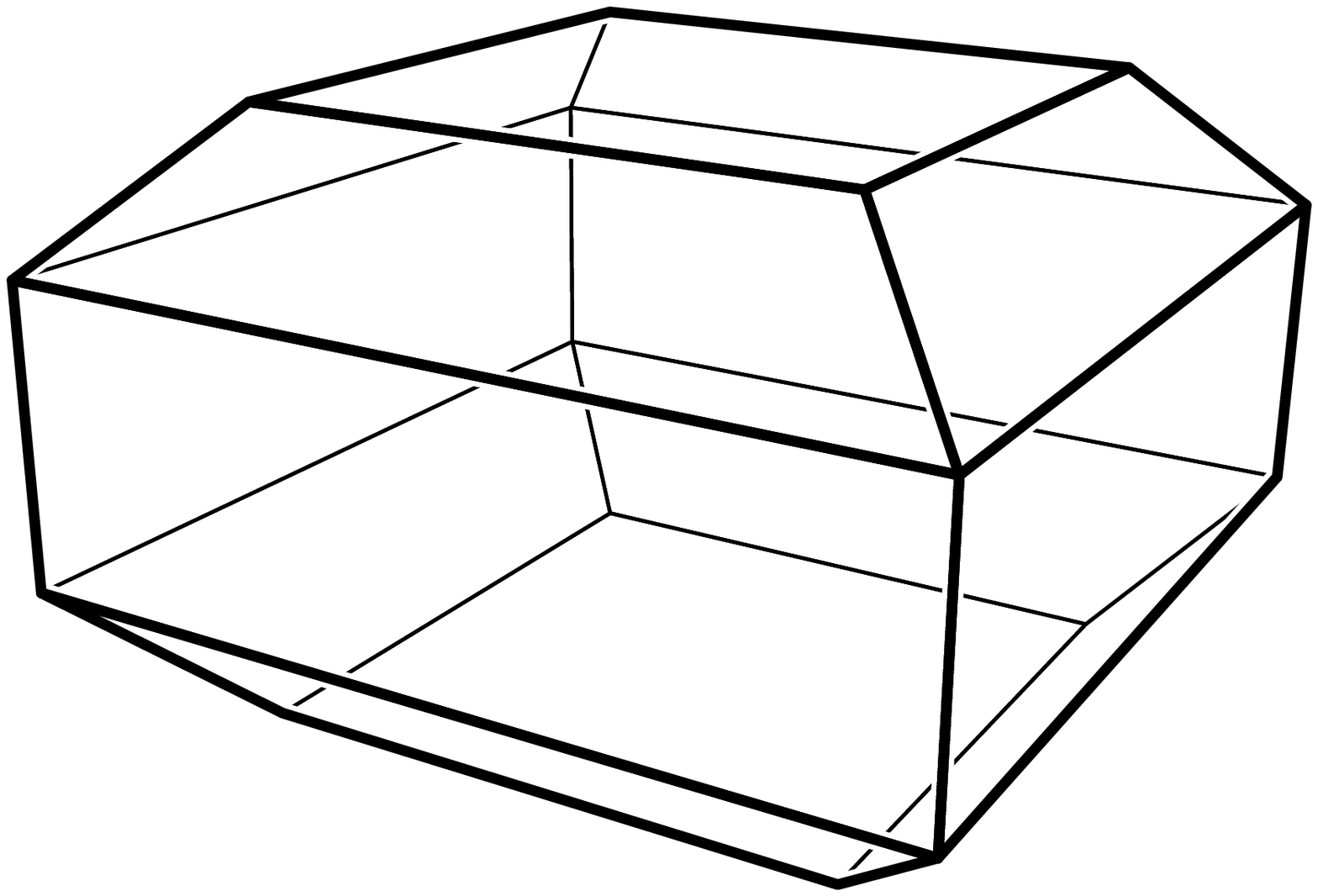,width=6cm}
    \end{minipage}
    \begin{minipage}[c]{6cm}
      \hspace*{1cm}\epsfig{file=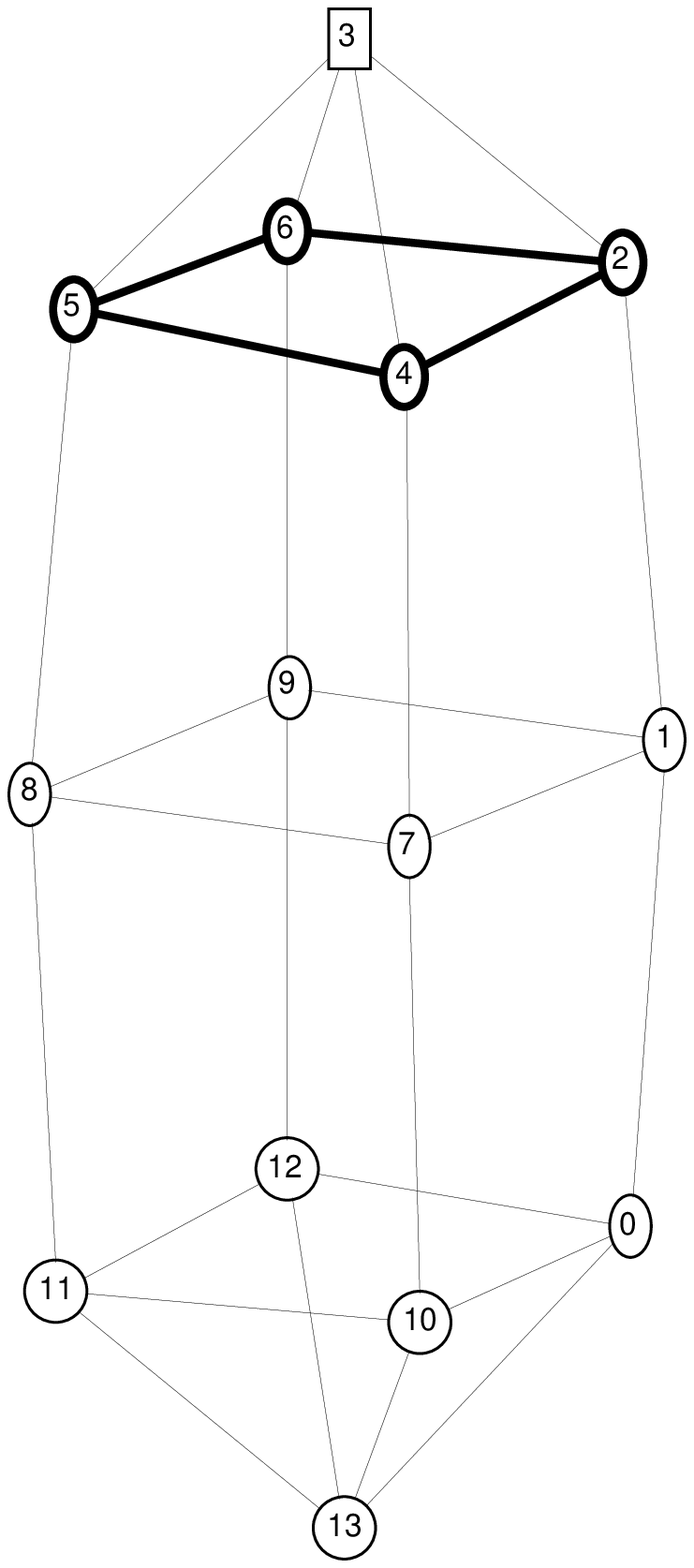,height=6cm}
    \end{minipage}
    \caption[Capped and not capped.]{A capped cubical polytope and its dual graph.}
    \label{fig:dual-graph}
  \end{center}
\end{figure}

By definition every facet of a cubical $d$-polytope is a $(d-1)$-cube.  In particular, every facet
contains $2(d-1)$ ridges, i.e.\ the dual graph of a cubical $d$-polytope is $2(d-1)$-regular.  A
vertex figure of any vertex in a dual-to-cubical $d$-polytope is a $(d-1)$-dimensional cross
polytope.  In view of Theorem~\ref{thm:MainTheorem} in order to prove the conjecture it is
sufficient to identify how the neighbors of any facet in the dual graph can be related to the
vertices of a cross polytope.  The main problem is that the induced subgraph among the neighbors in
the dual graph is not the graph of a cross polytope, in general.  Usually, there are many edges
missing.  For instance, in the example of a cubical $3$-polytope in Figure~\ref{fig:dual-graph} the
induced subgraphs of the facets numbered $1$, $7$, $8$, and $9$, are totally disconnected.

We will show that Conjecture~\ref{Con} holds for a very special class of cubical polytopes: A
cubical polytope~$P$ is called {\em capped\/} over a cubical polytope~$Q$ if there is a
combinatorial cube~$C$ such that $P=Q\cup C$ and $Q\cap C$ is a facet of~$Q$.  The unique facet
of~$P$ which does not contain any vertex of~$Q$ is called the {\em cap\/} of~$P$ with respect
to~$Q$.  A cubical polytope is called {\em capped\/} if it is obtained from iterated capping
starting with a combinatorial cube; see Blind and Blind~\cite{970.45465}.

\begin{figure}[htbp]
  \begin{center}
    \begin{minipage}[c]{6cm}
      \epsfig{file=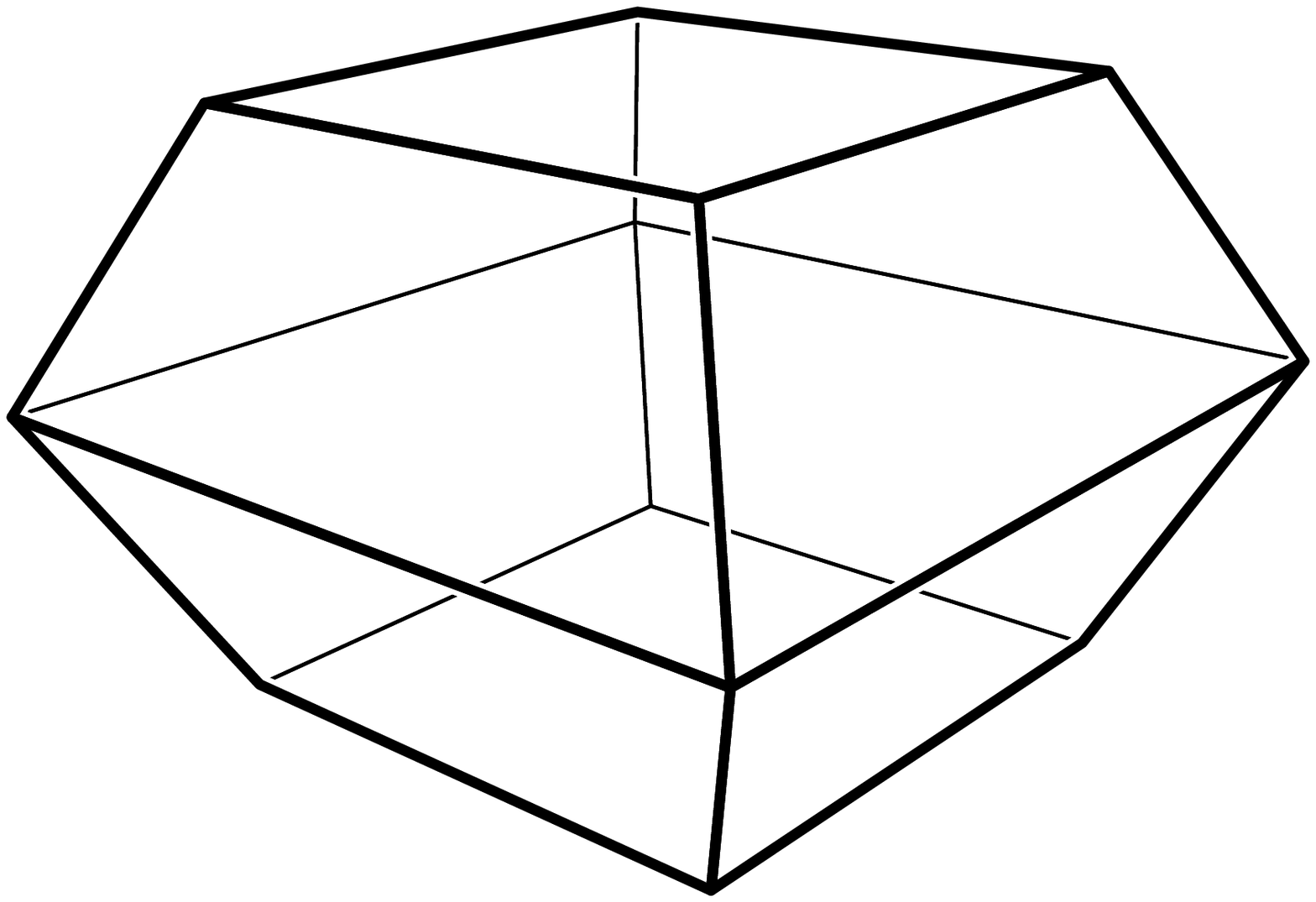,width=6cm}
    \end{minipage}
    \begin{minipage}[c]{6cm}
      \epsfig{file=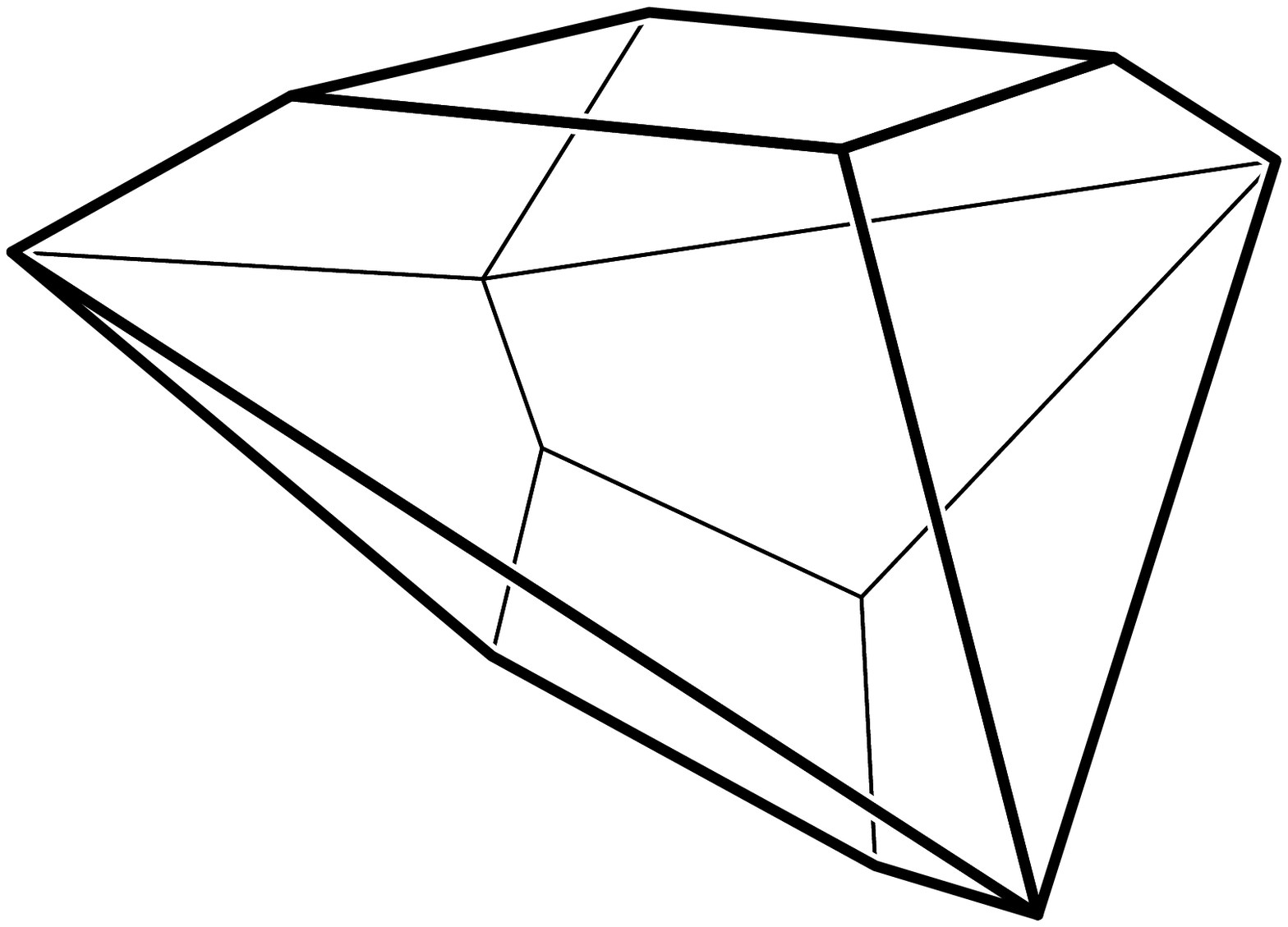,width=6cm}
    \end{minipage}
    \caption[Capped and not capped.]{A capped and a not capped cubical polytope which are
      combinatorially equivalent.  The polytope on the right is not capped because the quadrangle
      consisting of the $4$~vertices around the waist is not planar.}
    \label{fig:capped}
  \end{center}
\end{figure}

A few remarks on capped cubical polytopes.  The property of being capped is {\em not\/} a
combinatorial property as can be seen from Figure~\ref{fig:capped}.  For an example of a cubical
polytope which is not combinatorially equivalent to a capped polytope, see the zonotope in
Figure~\ref{fig:zonotope}.

Before we will study dual graphs of cubical polytopes it is helpful to collect a few results about
the {\em graphs\/} of cubical polytopes (and even more general graphs). 

A pure polytopal complex is called {\em cubical} if all its facets are combinatorial cubes.  A {\em
  cubical sphere\/} is a cubical complex which is homeomorphic to a sphere.  We want to introduce
the notion of {\em constructibility\/} of a pure cubical complex inductively.  Every cube is {\em
  constructible}.  If $A$ and $B$ are constructible $d$-dimensional cubical complexes, $A\cup B$ is
a cubical complex, and $A\cap B$ is pure $(d-1)$-dimensional constructible, then $A\cup B$ is {\em
  constructible}.  The {\em graph\/} of a complex is its $1$-skeleton.  Note that a constructible
complex (and thus also its graph) are necessarily connected, as a straightforward inductive argument
shows.

The following observation is due to G\"unter M. Ziegler.

\begin{prp}\label{prp:no-odd-cycles}
  The graph of a constructible cubical complex of dimension at least~$2$ does not contain any odd
  cycles.  In particular, its graph is bipartite.
\end{prp}

\begin{proof}
  No cube contains any odd cycle.
  
  Assume there is a constructible cubical complex~$C$ which contains an odd cycle~$\gamma$ and which
  is not a cube.  Then there are constructible cubical complexes~$A,B$ such that $A\cup B=C$.  By an
  induction on the number of construction steps necessary to build up~$C$ we can assume that
  $\gamma$ is contained in neither~$A$ nor~$B$.  That is, $\gamma$ passes through the
  intersection~$A\cap B$, which is connected, because of our assumption on the dimension.  It is
  conceivable that $\gamma$ enters and leaves~$A\cap B$ several times.  Pick two vertices~$x,y$
  on~$\gamma$ in~$A\cap B$ such that neither half of~$\gamma$, obtained from cutting at~$x$ and~$y$,
  is contained in~$A\cap B$.  If there is only one vertex in the intersection, then either $A$ or
  $B$ contains an odd cycle and we are done.  Choose a path~$\pi$ in $A\cap B$ between $x$ and~$y$.
  Depending on the parity of the length of~$\pi$, combining $\pi$ with either half of~$\gamma$
  yields an odd cycle~$\gamma'$.  Observe that $\gamma'$ enters and leaves the intersection one
  times less than~$\gamma$.  An obvious induction now gives the result.
\end{proof}

The boundary of a triangle is a $1$-dimensional constructible cubical complex.  This shows that the
assumption concerning the dimension is necessary.

By a theorem of Bruggesser and Mani~\cite{251.52013} the boundary of a polytope is known to be
shellable and thus constructible.  In particular, the boundary of a cubical polytope is a
constructible cubical complex.  So the above result implies that the graph of a cubical
$d$-polytope, where $d\ge 3$, does not contain any triangle.

Two disjoint facets of a cube are called {\em opposite}.

\begin{lem}\label{lem:3-in-a-line-opp}
  Let $A$, $B$, $C$ be facets in a cubical $d$-polytope~$P$, where $d\ge3$, such that $A\cap B$ and
  $B\cap C$ are opposite ridges in~$B$.
  
  Then $A\cap C$ is not a ridge.
\end{lem}

\begin{figure}
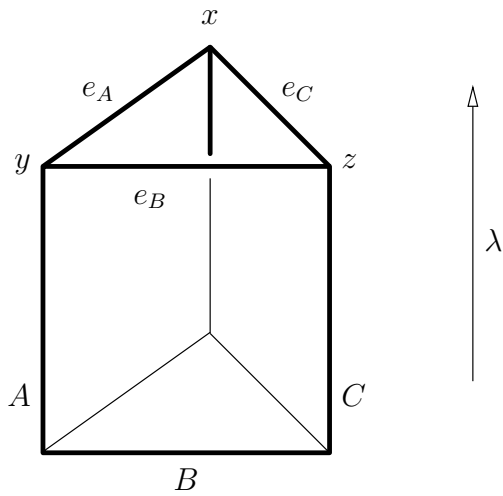

  \begin{center}
    \input 3-facets.pstex_t
    \caption{Three cubical facets and a linear function.}
    \label{fig:3-facets}
  \end{center}
\end{figure}

\begin{proof}
  Assume $A\cap C$ {\em is\/} a ridge, cf.\ Figure~\ref{fig:3-facets}.
  
  Choose an edge~$e_A$ in~$A$ which is not contained in~$B$ or~$C$.  Let $(x,y)=e_A$ with $x\in C$
  and $y\in B$.  There is a unique edge~$e_B$ through~$y$ which is contained in~$B$, but not
  contained in~$A$.  Now choose a linear objective function~$\lambda$ such that:
  \begin{enumerate}
  \item the edge $e_B$ is the unique maximal face of~$B$ with respect to~$\lambda$, and
  \item the function $\lambda$ attains the same value on~$e_B$ and at~$x$.
  \end{enumerate}
  
  Let $m_A$ be the maximal face of~$A$ with respect to~$\lambda$.  Clearly, $m_A$ contains the
  edge~$e_A$.  If $m_A$ were strictly greater than~$e_A$, then the intersection $m_A\cap B$ would be
  strictly greater than~$y$.  This would contradict the maximality of~$e_B$.  We have $m_A=e_A$.
  
  Denote by~$z$ the vertex of~$e_B$ which is contained in~$C$.  Let $m_C$ be the maximal face of~$C$
  with respect to~$\lambda$.  Now $m_C$ contains the points $x$ and~$z$.  As in the argument above,
  $m_C\supsetneqq\{x,z\}$ contradicts the maximality of~$e_B$ in~$B$.  We conclude that
  $m_C=\{x,z\}$ is an edge.  The vertices $x$, $y$, $z$ form a triangle in the graph of the cubical
  polytope~$P$.  Due to Proposition~\ref{prp:no-odd-cycles} we arrive at the final contradiction.
\end{proof}

\begin{figure}
  \begin{center}
    \input{Moebius.pstex_t}
    \caption{A cubical complex homeomorphic to the M\"obius strip.}
    \label{fig:Moebius}
  \end{center}
\end{figure}
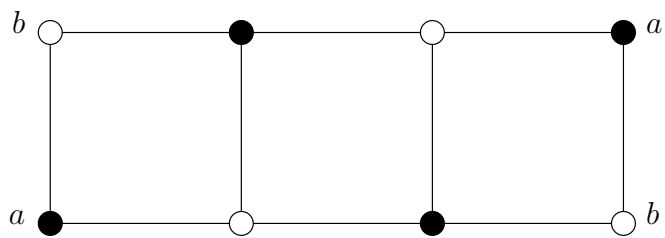

Note that there {\em are\/} (abstract) cubical complexes with a bipartite graph which do violate the
conclusion of the preceding lemma: Take a M\"obius strip built out of three quadrangles, cf.\
Figure~\ref{fig:Moebius}.

\begin{figure}
  \begin{center}
    \epsfig{file=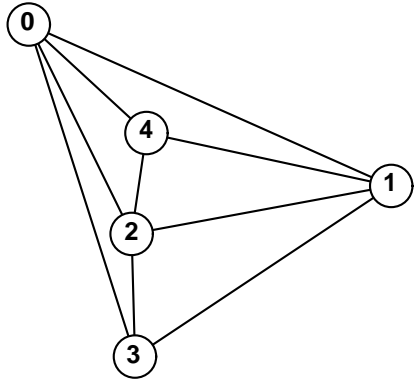}
    \caption{The graph of a $3$-polytope with edge~$(1,2)$ induced among the
      neighbors of~$0$.}
    \label{fig:diagonal}
  \end{center}
\end{figure}

Dualizing Lemma~\ref{lem:3-in-a-line-opp} yields that, if $e$ is an edge between two neighbors of a
vertex~$v$ in a dually cubical polytope, then $e$ and $v$ form a triangular face.  For general
polytopes it may happen that the induced subgraph on the neighbors of a given vertex~$v$ has edges
which are not contained in any proper face through~$v$, see the example in
Figure~\ref{fig:diagonal}.

A further consequence of the above lemma is that the induced subgraph among the neighbors in the
dual graph of a cubical $d$-polytope is a subgraph of a complete graph on $2(d-1)$~vertices minus a
perfect matching.  This means that finding the edge labeling of the vertex figure of a facet in the
dual graph is particularly easy for a cubical polytope.

\begin{lem}\label{lem:opposite}
  Let $P$ be a cubical $d$-polytope, $F$ one of its facets, $\Omega$ the set of neighbors of~$F$
  in~$\Delta_P$.
  
  Then the set
  $$\SetOf{\{N_1,\ldots,N_{d-1}\}}{N_i\in\Omega,\ \mbox{\rm any two facets $N_i$ and $N_j$ intersect
      non-trivially in~$P$}}$$
  is the edge labeling of the vertex figure of~$F$ in~$P^\dual$.
\end{lem}

Differently phrased, Lemma~\ref{lem:opposite} says that it suffices to recognize the antipodal pairs
among the neighbors of a facet in the dual graph.  A similar result holds for uniform oriented
matroids, that is, oriented matroids which generalize cubical zonotopes; cf.\ Cordovil, Fukuda, and
Guedes de Oliveira~\cite[Theorem~3.1]{CordovilFukudaOliveira}.

\begin{lem}\label{lem:cross-polytope-graph}
  Let $P$ be capped over the cubical polytope~$Q$ with cap~$F$.
  
  Then the induced subgraph~$\Omega$ among the neighbors of~$F$ in~$\Delta_P$ is a complete graph minus
  a perfect matching.  Moreover, $\Delta_Q$ is obtained from $\Delta_P$ by
  contracting~$\Omega\cup\{F\}$ to a single node.  The edge labeled vertex figures of~$Q^\dual$
  determines the edge labeled vertex figures of~$P^\dual$ and vice versa.
\end{lem}

In Figure~\ref{fig:dual-graph} the facet numbered~$3$ is a cap.  The induced subgraph among its
neighbors is the graph of a quadrangle which is a $2$-dimensional cross polytope.

If we already know that the cubical polytope is capped, then we also obtain the converse of
Lemma~\ref{lem:cross-polytope-graph}.

\begin{lem}\label{lem:capped}
  Let $P$ be a capped cubical $d$-polytope.
  
  If $F$ is a facet with the property that the induced subgraph among the neighbors of~$F$
  in~$\Delta_P$ is a complete graph minus a perfect matching, then $F$ is a cap.
\end{lem}

\begin{proof}
  Consider the set~$\cR$ of all the ridges contained in neighbors of~$F$ in~$\Delta_P$, but which
  are disjoint from~$F$.  By assumption the sublattice of the face lattice of~$P$ generated by~$\cR$
  is isomorphic to the face lattice of a $(d-1)$-cube.  We have to show that the ridges of~$F$ are
  the facets of a convex (combinatorial) cube, i.e.\ we have to show that the ridges in~$\cR$ are
  contained in a hyperplane.

  A capped cubical $d$-polytope can also be seen as a $d$-dimensional (constructible) cubical
  complex.  The facets of this complex are the (combinatorial) cubes which are attached one by one.
  Each facet of a capped cubical polytope is contained in a unique facet of the complex, i.e.\ a
  unique $d$-cube.  Say, $F$ is a facet of the $d$-cube~$C$.  Then the ridges in~$\cR$ are contained 
  in the unique facet of~$C$ which is opposite to~$F$.
\end{proof}

If we have a capped cubical polytope~$P$, then Lemma~\ref{lem:capped} allows us to detect a cap by
looking at the dual graph.  Applying Lemma~\ref{lem:cross-polytope-graph} then allows us to remove
the cap, thus yielding a capped cubical polytope with fewer facets.  Iterating this procedure we
gather the edge labeled vertex figures of~$P$.  Now Theorem~\ref{thm:MainTheorem} gives the desired
result.

\begin{thm}\label{thm:capped}
  Every capped cubical polytope can be reconstructed from its dual graph.
\end{thm}

Note the property in Lemma~\ref{lem:capped} does not characterize a capping.  If the induced
subgraph among the neighbors of some facet~$F$ happens to be complete minus a perfect matching, then
the ridges contained in the neighbors of~$F$ which are disjoint from~$F$ are combinatorially
isomorphic to the boundary of a cube, but they are not necessarily contained in a hyperplane.  An
example for this situation can be seen in Figure~\ref{fig:capped} (right).

\section*{Acknowledgements}

I am indebted to Manoj K.~Chari, Volker Kaibel and G\"unter M.~Ziegler for various helpful
suggestions.  The polytope pictures have been produced with Geomview~\cite{geomview} and
Graphlet~\cite{graphlet} via {\tt polymake}~\cite{polymake}.\nocite{DMV-Sem-1997}

\bibliographystyle{amsplain}
\bibliography{mic,math,top,poly,geo,soft,om}

\noindent{\small%
Michael Joswig\\
Fachbereich Mathematik, Sekr.\ 7--1\\
Technische Universit\"at Berlin\\
Stra\ss{}e des 17.~Juni~136\\
10623 Berlin, Germany\\
{\tt joswig@math.tu-berlin.de}}

\end{document}

%% file: 3-facets.pstex_t
\begin{picture}(0,0)%
\special{psfile=3-facets.pstex}%
\end{picture}%
\setlength{\unitlength}{3947sp}%
\begingroup\makeatletter\ifx\SetFigFont\undefined
\def\x#1#2#3#4#5#6#7\relax{\def\x{#1#2#3#4#5#6}}%
\expandafter\x\fmtname xxxxxx\relax \def\y{splain}%
\ifx\x\y   
\gdef\SetFigFont#1#2#3{%
  \ifnum #1<17\tiny\else \ifnum #1<20\small\else
  \ifnum #1<24\normalsize\else \ifnum #1<29\large\else
  \ifnum #1<34\Large\else \ifnum #1<41\LARGE\else
     \huge\fi\fi\fi\fi\fi\fi
  \csname #3\endcsname}%
\else
\gdef\SetFigFont#1#2#3{\begingroup
  \count@#1\relax \ifnum 25<\count@\count@25\fi
  \def\x{\endgroup\@setsize\SetFigFont{#2pt}}%
  \expandafter\x
    \csname \romannumeral\the\count@ pt\expandafter\endcsname
    \csname @\romannumeral\the\count@ pt\endcsname
  \csname #3\endcsname}%
\fi
\fi\endgroup
\begin{picture}(3174,3090)(3202,-4336)
\put(3526,-2311){\makebox(0,0)[rb]{\smash{\SetFigFont{12}{14.4}{rm}$y$}}}
\put(4276,-2536){\makebox(0,0)[b]{\smash{\SetFigFont{12}{14.4}{rm}$e_B$}}}
\put(4651,-1411){\makebox(0,0)[b]{\smash{\SetFigFont{12}{14.4}{rm}$x$}}}
\put(5476,-2311){\makebox(0,0)[lb]{\smash{\SetFigFont{12}{14.4}{rm}$z$}}}
\put(4051,-1861){\makebox(0,0)[rb]{\smash{\SetFigFont{12}{14.4}{rm}$e_A$}}}
\put(5101,-1861){\makebox(0,0)[lb]{\smash{\SetFigFont{12}{14.4}{rm}$e_C$}}}
\put(6376,-2836){\makebox(0,0)[lb]{\smash{\SetFigFont{12}{14.4}{rm}$\lambda$}}}
\put(4501,-4336){\makebox(0,0)[b]{\smash{\SetFigFont{12}{14.4}{rm}$B$}}}
\put(3526,-3811){\makebox(0,0)[rb]{\smash{\SetFigFont{12}{14.4}{rm}$A$}}}
\put(5476,-3811){\makebox(0,0)[lb]{\smash{\SetFigFont{12}{14.4}{rm}$C$}}}
\end{picture}

%% file: Moebius.pstex_t
\begin{picture}(0,0)%
\special{psfile=Moebius.pstex}%
\end{picture}%
\setlength{\unitlength}{3947sp}%
\begingroup\makeatletter\ifx\SetFigFont\undefined
\def\x#1#2#3#4#5#6#7\relax{\def\x{#1#2#3#4#5#6}}%
\expandafter\x\fmtname xxxxxx\relax \def\y{splain}%
\ifx\x\y   
\gdef\SetFigFont#1#2#3{%
  \ifnum #1<17\tiny\else \ifnum #1<20\small\else
  \ifnum #1<24\normalsize\else \ifnum #1<29\large\else
  \ifnum #1<34\Large\else \ifnum #1<41\LARGE\else
     \huge\fi\fi\fi\fi\fi\fi
  \csname #3\endcsname}%
\else
\gdef\SetFigFont#1#2#3{\begingroup
  \count@#1\relax \ifnum 25<\count@\count@25\fi
  \def\x{\endgroup\@setsize\SetFigFont{#2pt}}%
  \expandafter\x
    \csname \romannumeral\the\count@ pt\expandafter\endcsname
    \csname @\romannumeral\the\count@ pt\endcsname
  \csname #3\endcsname}%
\fi
\fi\endgroup
\begin{picture}(4177,1448)(1974,-3894)
\put(2251,-3811){\makebox(0,0)[rb]{\smash{\SetFigFont{12}{14.4}{rm}$a$}}}
\put(6151,-3811){\makebox(0,0)[lb]{\smash{\SetFigFont{12}{14.4}{rm}$b$}}}
\put(6151,-2611){\makebox(0,0)[lb]{\smash{\SetFigFont{12}{14.4}{rm}$a$}}}
\put(2251,-2611){\makebox(0,0)[rb]{\smash{\SetFigFont{12}{14.4}{rm}$b$}}}
\end{picture}